\documentclass[reqno, 12pt]{amsart}

\def\C{{\mathbb C}}

\def\Z{{\mathbb Z}}
\def\Zp{{\mathbb Z}_p}

\def\A{{\bf A}}
\def\K{{\bf K}}
\def\k{{\bf k}}
\def\sep{{{\rm sep}}}

\def\Fp{{{\mathbb F}_p}}
\def\F{{\mathbb F}}
\def\Fq{{{\mathbb F}_q}}
\def\Fr{{{\mathbb F}_r}}

\def\isom{\cong}
\def\dis{\displaystyle}
\def\scr{\scriptstyle}

\newtheorem{lemma}{Lemma}
\newtheorem{cor}{Corollary}
\newtheorem{prop}{Proposition}

\newtheorem{theorem}{Theorem}

\theoremstyle{definition}
\newtheorem{defn}{Definition}

\newtheorem{questions}{Question}

\theoremstyle{remark}
	
\newtheorem{rems}{Remarks}	

\newtheorem{examples}{Example}

\begin{document}

\title[Separability, multi-valued operators, and zeroes of $L$-functions]{Separability, multi-valued operators,\\ and zeroes of $L$-functions}
\author{David Goss}
\thanks{This research was supported by NSA grant MDA 904-97-1-0048}
\address{Department of Mathematics\\ The Ohio State University\\ 231 W.
$18^{\text{th}}$ Ave. \\ Columbus, Ohio 43210}
\email{goss@math.ohio-state.edu}
\date{Autumn, 1997}

\begin{abstract}Let $\k$ be a global function field in $1$-variable over
a finite extension of $\Fp$, $p$ prime, $\infty$ a fixed place of
$\k$, and $\A$ the ring of functions of $\k$ regular outside of
$\infty$. Let $E$ be a Drinfeld module or $T$-module. Then,
as in \cite{go1}, one can construct associated characteristic $p$
$L$-functions based on the classical model of abelian
varieties {\it once} certain auxiliary choices are made. Our purpose
in this paper is to show how the well-known concept of
``maximal separable (over the completion $\k_\infty$) subfield''
allows one to construct from such $L$-functions 
certain separable extensions which are independent of these choices.
These fields will then depend only on the isogeny class
of the original $T$-module or Drinfeld module and $y\in \Zp$,
and should presumably
be describable in these terms. Moreover, they give a very useful framework
in which to view the ``Riemann hypothesis'' evidence of
\cite{w1}, \cite{dv1}, \cite{sh1}.
We also establish that an element which is
{\it separably} algebraic over $\k_\infty$ can be realized as
a ``multi-valued operator'' on general $T$-modules.
This is very similar to realizing $1/2$ as the 
multi-valued operator $x\mapsto \sqrt{x}$ on $\C^\ast$. Simple examples
show that this result is false for non-separable elements. This result may
eventually allow a ``two $T$'s'' interpretation of the above extensions in
terms of multi-valued operators on $E$ and certain tensor twists.

\end{abstract}

\maketitle


\section{Introduction}\label{intro}
Let $\Fr$ be the finite field with $r=p^m$ elements and let 
$\mathcal X$ be a smooth projective geometrically connected curve over
$\Fr$. Let $\infty\in \mathcal X$ be a fixed closed point and let
${\rm Spec}(\A):={\mathcal X}\setminus \infty$. 
It is well known that $\A$ is a Dedekind domain
with a finite class group and unit group equal to $\Fr^{\ast}$. The domain
$\A$ is also the ground ring for the theory of Drinfeld modules and
$T$-modules --- in this paper we shall loosely refer to Drinfeld modules or $T$-modules
as ``motives.'' To such motives one can associate
characteristic $p$ valued $L$-series via the familiar use of Euler factors.
These concepts were discussed in detail in \cite{go1}; in 
particular, attention was given in \cite{go1} to some numerical
evidence \cite{w1}, \cite{dv1} (and now \cite{sh1}) which seems
to be an indication of the correct
``Riemann hypothesis'' in the theory. In Section 8.24 of \cite{go1}
a rather ad hoc
exposition of the available evidence was presented.

Another major theme of \cite{go1} was the ``two $T$'s'' which we now
briefly explain. Let $\k$ be the fraction field of $\mathcal X$ (= the quotient
field of $\A$) and let $F/\k$ be a finite extension. Let $\psi$ be a Drinfeld
module over $F$. Then via $\psi$ the elements of $\A$ play {\em two}
distinct roles in the theory: sitting inside $F$, they are {\em scalars}, while
inside $\A$ they are {\it operators}. The ``two $T$'s'' is just the
idea that one must constantly keep these two concepts distinct. Thus one works
with algebras like $\A\otimes_{\Fr}\A$ where one copy of $\A$ represents
scalars and the other represents operators. (In the paper we will denote
the copy representing scalars with non-bold symbols to avoid any confusion.)

In any case, simple
as it may seem, the two $T$'s is a very useful
principle for grouping the various elements of function field arithmetic. Indeed,
in the universe of $T$ as scalar one finds exponential functions, periods,
modular forms, factorials,
gamma functions, etc. In the universe of $T$ as operator one
has  the theory of the zeta function of Drinfeld modules and
$T$-modules over finite field,
global $L$-functions of Drinfeld modules or $T$-modules, zeta functions, etc.
Indeed, it is precisely the $\A$-action on the motive that allows one to
form characteristic polynomials of the Frobenius morphism of a $T$-module
over a finite field and these characteristic polynomials (which 
live the in operator universe)
are precisely the Euler factors of the associated global
$L$-series. The two $T$'s also seems to
suggest viewing the zeroes of these $L$-functions as being operators
of some sort on the original motive.

The reader may well wonder about an analog of the two $T$'s in classical
algebraic number theory. However, 
as pointed out in Remark 9.9.13 of \cite{go1} this is
not possible simply because $\Z\otimes_{\Z}\Z$ {\it is} $\Z$. Thus complex
numbers must play the two roles at the same time; a sort of ``wave-particle
duality'' in number theory.

Now let $f(x)$ be an irreducible polynomial with coefficients in
some field of finite characteristic. It has long been known that the classical
methods of Galois theory, i.e., field automorphisms, work best only
for those $f(x)$ where $f(x)$ is {\it separable}. It is the purpose of
this note to begin the study of $L$-zeroes from the point of view
of the two $T$'s via the concept of separability. More precisely, we have
two major goals in this paper. The first will be to use the
concept of separability to present a general structure in which to view
the above mentioned ``RH''-evidence. This structure arises very simply
in the following manner: in order to even define $L$-functions of
motives (which will be recalled in Section \ref{lfunction}) there are a
couple of arbitrary choices that must be made --- one must 
first choose a sign function,
sgn, and then a positive uniformizer $\pi$ in the completion of $\k$ at
$\infty$. Once these two choices are made, one can ``exponentiate
ideals'' (Definition \ref{exponen}) and thus define $L$-functions on the
space $S_\infty:={\mathbf C}_\infty^\ast \times \Zp$, where 
${\mathbf C}_\infty$ is the completion of a fixed algebraic closure of
$\k_\infty$. Thus the zeroes of our $L$-functions at a given $y\in
\Zp$ (which
come from the first variable of $s=(x,y)\in S_\infty$)
actually depend on two auxiliary parameters: sgn and $\pi$.
However,
by simply passing to the {\it maximal separable (over $\k_\infty$) subfields}
of the fields generated by the zeroes for each $y$, we will show how we
already obtain an invariant of these choices. We call this separable
field $Z(L,y)$. There is then an ``obvious'' choice for this field
which arises simply from the coefficients of the $L$-series; this field
is called $C(L,y)$ and always lies in $Z(L,y)$. 
The RH-evidence seems to imply that we should
expect $C(L,y)=Z(L,y)$ (but see Remarks \ref{coefffield} for the
case of complex multiplication). In any case the fields $Z(L,y)$ depend only
on the isogeny class of the original motive and $y\in \Zp$ and should somehow
have a canonical description in these terms. 

The reader may also wonder about the total field generated by the roots
at a given $y\in \Zp$. However, it is very easy to see that a given
extension of $\k_\infty$ is determined by its maximal separable subfield {\it and}
the inseparability degree. Thus the crucial work is to find $Z(L,y)$.

Our second goal is to show how the concept of separability allows one to
extend the definition of the Drinfeld module or $T$-module to certain
``multi-valued operators'' as implied by the two $T$'s (and as mentioned
above). Let $\alpha$ be a complex number; then the
operator $x\mapsto x^\alpha=e^{\alpha \log x}$ on ${\mathbf C}^\ast$ is well
known to be multi-valued in general. 
By using the logarithm and exponential of a
motive one can try to
 define multi-valued operators in characteristic $p$ in a very
similar fashion. In fact it turns out that it is impossible in
general (i.e., for arbitrary $T$-modules) to realize a given
$\alpha\in {\mathbf C}_\infty$ as a multi-valued operator as simple examples
show. However, this procedure will be shown to {\it always} work
when $\alpha$ is
separably algebraic over $\k_\infty$ (Theorem \ref{mainth}; see also
the remark before Question \ref{multicompl} for the case of c.m.).
It thus remains to somehow connect the two uses of separability in this
paper; i.e., to canonically describe the fields $Z(L,y)$ and
$C(L,y)$ etc. in terms of
multi-valued operators. As of this writing, there are no obvious clues
on how to do this and we pose several questions along these lines
(e.g. Questions \ref{bigbasic} and \ref{multicompl}).

Section \ref{lfunction} reviews the definitions of $L$-series and
defines the fields $Z(L,y)$, $C(L,y)$, etc., (Definitions \ref{zfields} and
\ref{defcl}). In Remarks \ref{nofe} we will also show how the two $T$'s gives
a satisfactory explanation of why the characteristic $p$ $L$-functions do not
satisfy functional equations of the classical sort.
Section \ref{drinfeld} then discusses multi-valued operators
and Drinfeld modules where the desired realization as multi-valued
operators is very easy. Section \ref{tmodule} discusses the general theory
of multi-valued operators for arbitrary $T$-modules. We will show that there
is a very close connection between realizing elements as multi-valued
operators and the theory of hyper-derivatives. Finally Section
\ref{vadic} discusses what little is known for the $v$-adic theory
of $L$-functions of Drinfeld and $T$-modules, $v\in {\rm Spec}(\A)$.
In particular we also report on some interesting computer work
on certain {\it global} fields (i.e., finite
extensions of $\k$) which arise in the theory.

It is my pleasure to thank the National Security Agency for its support
of this work and to acknowledge very useful communication with G. Anderson,
T. Satoh and D. Wan. I also thank J. Roberts for his computer work which is
reported on in Section \ref{vadic} and J. Zhao for pointing out a number
of misprints in an earlier version of this paper.


\section{Separability and zeroes of $L$-functions}\label{lfunction}
In Section 8 of \cite{go1} the theory of characteristic $p$ valued 
$L$-functions was presented, and, in particular, some evidence for
a ``Riemann hypothesis'' for such functions was given (Section 8.24). 
It is the goal of this section to recast this evidence from the viewpoint of
``separable subfields of zero fields'' and
the ``two $T$'s'' which we believe is a more natural approach than the ad hoc
one taken in \cite{go1}. In particular, we will formulate 
here some general problems
on the zeroes of arbitrary $L$-functions; the solution to these problems
may also explain the above mentioned evidence.

We begin by reviewing the definition of our $L$-series where the basic
notation of rings will be exactly as in \cite{go1}.
We let $\A$ be arbitrary, as in the introduction, and set $\k$ to
be its quotient field. Set $d_{\infty}=\deg \infty$ where the degree is taken over
$\Fr$ and let $\K=\k_\infty$. Let ${\mathbf C}_\infty$ be the completion of
a fixed algebraic closure $\overline \K$  of $\K$ equipped with its canonical metric, and let
$\overline{\k}\subset {\mathbf C}_\infty$ be the algebraic closure of $\k$. Let
$\k^{\rm sep}\subset \overline{\k}$ be the separable closure of $\k$ and let
$\K^{\rm sep}\subset {\bf C}_{\infty}$ be the separable closure of $\K$.
Let $\F_\infty\subset \K$ be the field of constants.
A {\it sign function}, sgn, is a homomorphism ${\rm sgn}\colon \K^\ast \to
\F_\infty^\ast$ which is the identity on $\F_\infty^\ast$. An element
$\alpha\in \K^\ast$ is said to be {\it positive} or {\it monic} if and only
if ${\rm sgn}(\alpha)=1$.

Now fix a positive uniformizer $\pi\in \K^\ast$ and let $\alpha$ be positive.
We therefore have
$$\alpha =\pi^j \langle \alpha \rangle\,,$$
where $j=v_\infty(\alpha)$ and $\langle \alpha \rangle=\langle \alpha
\rangle_{\pi}$ depends on $\pi$ and belongs to the
$1$-units, $U_1$, in $\K$. The element
$\langle \alpha \rangle$ is called the {\it $1$-unit
part of $\alpha$ with respect to $\pi$}.
Clearly $\langle \alpha \beta \rangle=
\langle \alpha \rangle \langle \beta \rangle$.
Let $\mathcal I$ be the group of $\A$-fractional ideals and let
${\mathcal P}^+$ be the subgroup of principal and positively generated
ideals. It is quite easy to see that the finiteness of the class
group of $\A$ implies that ${\mathcal I}/{\mathcal P}^+$ is finite.
Let ${\mathfrak I}=({\mathfrak i})$ be in ${\mathcal P}^+$ where $\mathfrak i$
is positive. We then naturally define $\langle {\mathfrak I}\rangle:=
\langle {\mathfrak i} \rangle$ thus defining a homomorphism from
${\mathcal P}^+\to U_1.$ 

It is well known, and easy to see, that $U_1$ is a topological module over
the $p$-adic integers $\Zp$. Let $\hat{U_1}\subset {\mathbf C}_\infty$ be the
$1$-units of ${\mathbf C}_\infty$. It is clear that $\hat{U_1}$ also contains
unique $p$-th roots. Thus $\hat{U_1}$ is a $\mathbb{Q}_p$-vector space
and, consequently, is a uniquely divisible abelian group. As divisible abelian
groups are injective, the homomorphism $\langle~\rangle$ may be extended 
automatically to $\mathcal I$; as ${\mathcal I}/{\mathcal P}^+$ is finite,
this extension is {\it unique}. We will also denote this extended mapping
by $\langle ~\rangle \colon {\mathcal I}\to \hat{U_1}$.

\begin{defn} We set $S_\infty:={\mathbf C}^\ast_\infty\times \Zp$.
\end{defn}\noindent
The group $S_\infty$ has a topology in the obvious fashion; its group
operation is written additively. 

The next definition presents the fundamental notion of ``exponentiating
an ideal.''

\begin{defn}\label{exponen}
Let $\mathfrak I$ be a non-zero
fractional ideal of $\A$ and let $s=(x,y)\in S_\infty$. We set
$${\mathfrak I}^s:=x^{\deg \mathfrak I}\langle {\mathfrak I} \rangle^y\,,$$
where 
$$\langle {\mathfrak I}\rangle^y=\sum_{j=0}^\infty \binom{y}{j}(\langle {\mathfrak I}\rangle-1)^j\,.$$ 
\end{defn}

It is easy to check that ${\mathfrak I}\mapsto {\mathfrak I}^s$ 
has all the usual properties of
exponentiation $n\mapsto n^s$ for $n$ a positive integer and
$s$ a complex number. The group $S_{\infty}$ is thus seen to play the
role classically played by the complex numbers; that is, our $L$-functions are
naturally defined on an open subset of $S_{\infty}$. 

Let $\pi_\ast\in {\mathbf C}_\infty$ be a fixed
$d_\infty$-th root of $\pi$. Let $j\in \Z$. We set
$$s_j=(\pi_\ast^{-j},j)\in S_\infty\,.$$
When no confusion will result, we simply denote $s_j$ by ``$j$.''

\begin{prop}\label{is}
Let ${\mathfrak I}=({\mathfrak i})$ be a
principal ideal. Then ${\mathfrak I}^{s_j}={\mathfrak 
i}^j/{\rm sgn}({\mathfrak i})^j$.
\end{prop}

\begin{proof} This follows immediately from Corollary 8.2.7 in
\cite{go1}.\end{proof}

\begin{cor}\label{iposs}
Let ${\mathfrak I}=({\mathfrak i})$ be as above where $\mathfrak i$ is now
positive. Then ${\mathfrak I}^{s_j}={\mathfrak i}^j$.
\end{cor}

Next we need to discuss how the above definitions depend on the choice of
positive uniformizer. Let $_1\pi$ and $_2\pi$ be two positive uniformizers
in $\K$. Let $u:={_1\pi}/{_2\pi}$; thus $u\in U_1$. Let $\langle ?\rangle_i$
be the $1$-unit parts defined with respect to $_i\pi$ for $i=1,2$.

\begin{prop}\label{dependpi}
Let $\mathfrak I$ be a fractional ideal of $\A$.
Then
$$\langle {\mathfrak I}\rangle_1=\left( u^{1/d_\infty}\right)^{\deg \mathfrak I}\cdot
\langle {\mathfrak I}\rangle_2\,,$$
where we take the unique $1$-unit root of $u$.
\end{prop}

\begin{proof}
This is Proposition 8.2.15 of \cite{go1}.\end{proof}

\begin{rems}\label{totins}
Let $\alpha$ be any $1$-unit in $\K$ and let $t\in {\mathbb Q}_p$. If 
$t\in \Zp$, then $\alpha^t$ is computed by the binomial theorem
exactly as in Definition \ref{exponen}, and clearly also lies in
$\K$. If $t\not\in {\mathbb Z}_p$ then
there is certainly a power of $p$, say $p^m$, such that $t_1:=p^mt\in
\Zp$. Thus 
$$\alpha^t=(\alpha^{t_1})^{p^{-m}}$$
and is thus {\it totally inseparable} over $\K$. This applies to
$u^{1/d_\infty}$ where $u$ is given just above. It also applies to
the elements $\langle {\mathfrak I}\rangle$ for $\mathfrak I$ an
ideal of $\A$. Indeed, if $o$ is the order of ${\mathcal I}/{\mathcal P}^+$,
then ${\mathfrak I}^o=({\mathfrak i})$ where $\mathfrak i$ is positive; 
so $\langle {\mathfrak i} \rangle\in \K$. Further,
$$\langle {\mathfrak I}\rangle =\langle{\mathfrak i}\rangle^{1/o}$$
and is thus totally inseparable over $\K$. The same result is
now easily seen to
be true for $u^{y/d_\infty}=(u^{1/d_\infty})^y$ and $\langle {\mathfrak I}
\rangle^y$ for any $y\in \Zp$.
\end{rems}

Now let $_i\pi_\ast$ be chosen $d_\infty$-th roots of $_i\pi$, $i=1,2$.
Let $\mathfrak I$ be an ideal and let ${\mathfrak I}_i^{s_1}:={\mathfrak I}^{s_1}$ be defined
with respect to $_i\pi$, $i=1,2$; i.e., by using $_i\pi$ in the definition
of $s_1\in S_\infty$.

\begin{prop}\label{idependpi}
There exists a $d_\infty$-th root of unity $\zeta$ which is
independent of $\mathfrak I$, such that
$${\mathfrak I}_1^{s_1}=\zeta^{\deg I}\cdot {\mathfrak I}_2^{s_1}\,.$$
\end{prop}

\begin{proof}
This is Proposition 8.2.16 of \cite{go1}\end{proof}

Next we recall the definition of $L$-series of Drinfeld modules and
$T$-modules. As  stated in the introduction, a nonconstant element 
$T\in \A$ plays two distinct roles in the theory. In order to keep
these roles distinct, it is very convenient to introduce another
copy of the rings $\A$, $\k$, $\k_v$ ($v\in {\rm Spec}(\A)$), etc. --- to keep the two copies apart,
we will denote the second copy by non-bold letters. So the rings
$A$, $k$, $k_{v}$ ($v$ is identified with its isomorphic image in
${\rm Spec}(A)$), etc., will be rings of {\it scalars}. There is obviously a 
canonical isomorphism from a boldface ring to its non-bold copy; this
isomorphism is denoted $\theta:=\theta_{\rm con}$.  Thus the map,
$$\A \stackrel{\theta_{\rm con}}{\longrightarrow} A \hookrightarrow k$$
makes any extension of $k$ automatically into an $\A$-field. 
{\it We will adopt this ``two $T$'s'' set-up throughout this paper. Thus 
Drinfeld modules, etc., will be defined over non-bold rings of scalars inside
$C_{\infty}$ and will give
rise to algebraic actions of elements in the bold rings of
operators.} The reader should note that because we have adopted the
two $T$'s throughout this paper, our notation for $L$-series
will occasionally be different than that of \cite{go1} (e.g., Example
\ref{zeta}).

It is standard to denote $\theta_{\rm con}(T)$ be $\theta$. However, we
shall also adopt the useful convention (as in \cite{go1}) that if $\alpha
\in {\mathbf C}_{\infty}$ then $\overline{\alpha}:=\theta (\alpha)\in C_{\infty}$
etc. Thus, for instance, $\theta(T)=\theta=\overline{T}$.

Now let $F$ be a finite extension of $k$ and 
let ${\mathbf F}\subset \overline{\k}$ be its inverse image under
$\theta_{\rm con}.$ Let $\psi$ be a Drinfeld
module of rank $d$ defined over $F$. Let ${\mathcal O}_F\subset F$ be the
ring of $A$-integers. For almost all primes $\mathfrak P$ of $A$, 
one can reduce $\psi$ to a Drinfeld module $\psi^{\mathfrak P}$
of rank $d$ over the finite
field ${\mathcal O}_F/\mathfrak P$ (which is obviously still an $\A$-field).
Thus there is a Frobenius endomorphism $F_{\mathfrak P}$ of
$\psi^{\mathfrak P}$, and one sets
$$P_{\mathfrak P}(u):=\det (1-uF_{\mathfrak P}\mid T_v(\psi^{\mathfrak P}))\,$$
where $T_v$ is the $v$-adic Tate module of $\psi^{\mathfrak P}$ for
$v\in {\rm Spec}(\A)$ prime to $\mathfrak P$. One knows from work of Drinfeld that 
this polynomial has coefficients in $\A$ (i.e., the coefficients are
{\it operators}) which are
independent of the choice of $v$ (\cite{go1} 4.12.12.2), and that its
roots satisfy the local ``Riemann hypothesis'' (\cite{go1} 4.12.8.5).
One thus forms the $L$-series
$$L(\psi, s):=\prod_{\mathfrak P}P_{\mathfrak P}(N{\mathfrak P}^{-s})^{-1}\,$$
where the product is taken over all such $\mathfrak P$ as above,
 $s\in S_\infty$, and $N\mathfrak P$ is the ideal norm from ${\mathcal O}_F$ to $\A$. The local Riemann hypothesis implies that this product converges
on a ``half-plane'' of $S_\infty$ (\cite{go1}, 8.6.9.1). Of course one would also
like to also have factors at the bad primes of ${\mathcal O}_F$. As this is
not important to us in this work, we refer the reader to Subsection 8.6
of \cite{go1} --- in fact, very little is known here in general
and finding the right factors is an important problem.

The analytic properties of $L(\psi, s)$ are always handled in the following
fashion. One expands out the Euler product for $L(\psi,s)$ to obtain
a ``Dirichlet series'' of the form $\sum_{\mathfrak I} a_{\mathfrak I}
{\mathfrak I}^{-s}$. In turn, the sum $\sum a_{\mathfrak I}{\mathfrak I}^{-s}$
is rewritten as
$$\sum_{j=0}^\infty x^{-j}\left(\sum_{\deg {\mathfrak I}=j}a_{\mathfrak I}
\langle {\mathfrak I} \rangle^{-y}\right)\,.$$
The analytic properties of $L(\psi,s)$ are then determined by the 
convergence properties, etc., of this 
$1$-parameter family of power series in $x^{-1}$.

The procedure to define the $L$-series  of a general $T$-module $E$
with complex multiplication by
$\A$ (or ``$\A$-module'') is clear. One reduces $E$ at those primes where the
rank of the reduction remains the same, etc. The details have, as yet,
only been written down for the case $\A=\Fr[T]$ and we refer the
reader to \cite{go1} for this. In order to get good estimates on the
eigenvalues of the various Frobenius operators which is necessary
for $L$-series, we shall always
assume that our $T$-modules are pure in Anderson's sense
(see Definition 5.5.2 of \cite{go1}). For instance, tensor products of
Drinfeld modules are pure and the tensor product of two pure $T$-modules
is also pure, etc.

It is also clear, in analogy with classical theory, that the $L$-series of
a motive (i.e, a Drinfeld module or $\A$-module) defined over a field $F$ depends
only on the isogeny class of the motive over $F$.

One can more generally define $L$-series in the context of strictly
compatible families of $v$-adic representations of the Galois group
$F^{\rm sep}$ over $F$ (where $F^{\rm sep}\subset C_\infty$ is the
separable closure) in the standard
fashion of elliptic curves, etc., (\cite{go1} Definition 8.6.7).
The simplest such compatible family is obviously the trivial $1$-dimensional
character and the $L$-function of this family is denoted
$\zeta_{{\mathcal O}_F}(s)$. The $L$-series that are defined in this fashion are
said to be ``of arithmetic interest.''

\begin{examples}\label{zeta}
When $F=k$, we obtain the zeta function $\zeta_A(s)$ of $\A$. Thus, by
definition,
$$\zeta_{ A}(s)=\prod_{{\mathfrak P}\in\, {\rm Spec}(A)}(1-{\mathfrak P}^{-s})^{-1}=\sum_{\mathfrak I}I^{-s}\,,$$
where $\mathfrak I$ runs over the ideals of $A$ and where we
have identified isomorphic (under $\theta$) ideals in $A$ and $\A$. When $\A=\Fr[T]$, then
$$\zeta_A(s)=\sum_{n~\rm monic}n^{-s}\,.$$
(Note that in \cite{go1}, these functions were denoted $\zeta_\A(s)$ etc. To
be consistent with the two $T$'s we have dropped the bold-face in the 
definition...)
\end{examples}

\begin{rems}\label{zetaL}
1. Let $\A=\Fr[T]$ and let $C$ be the Carlitz module over $F$.
Anderson's result on tensor product representations (see, e.g. 5.7.3 of \cite{go1}) implies that for $n>0$
$$L(C^{\otimes n},s)=\zeta_{{\mathcal O}_F}(s-n)\,.$$
Thus, even in this case, one can not just study zeta-functions in the
context of Drinfeld modules alone, but {\it must}
pass to the category of general $T$-modules.\\
2. More generally, the same arguments imply that for any motive $M$
over $\Fr[T]$ we have 
$$L(M\otimes C^{\otimes n},s)=L(M,s-n)\,.$$
3. On the other hand, let $\A$ be arbitrary. Let $\psi$ be a sign normalized
rank one Drinfeld module (or ``Hayes-module'') for $\A$; these are the
natural generalization of the Carlitz module. Then $\psi$ is defined
only over the ring of integers ${\mathcal O}^+$ of the Hayes normalizing field
${H}^+\subset C_\infty$ of $k$. In particular $H^+$ is a certain finite
extension of the Hilbert class field $H$ of $k$ which is totally split
at infinity. As such, $\psi$ will be defined over $k$ only when $\A$
is a principal ideal domain. Therefore, Part 1 generalizes to $L(\psi,s)$ {\it only} for
those fields $F$ with
$H^+\subseteq  F$. In general, one can view
$\zeta_A(s)$ as being associated to the finite {\it set} of absolute
isomorphism classes of Hayes-modules (see, e.g., \S 8.19 of \cite{go1}).
\end{rems}

The reader should be aware that, by construction, all of the above functions 
$L(s)$ will
have the following basic property: there exists a fixed finite
extension $\K_{1}$ of $\K$ (depending on $L(s)$ of course) 
such that for every $y\in \Zp$, the coefficients of
the power series $L(x,y)$ will lie in $\K_{1}$. Let $y\in \Zp$.
Then standard results in 
non-Archimedean function theory now imply that the zeroes (in $x$)
of $L(x,y)$ belong to the algebraic closure $\overline{\K}$ of $\K$ inside
${\mathbf C}_\infty$. This simple observation is crucial in all of what
follows.

\begin{defn} \label{zerofield}
The smallest extension of $\K$ containing the zeroes of $L(x,y)$ is called
{\it the zero field of $L(s)$ at $y$}. If for a given $y\in \Zp$
the function $L(x,y)$ is
a non-zero constant, then we define $Z(L,y)=\K$.
As the zero field of $L(s)$
possibly depends on the sign function and $\pi$ as well
as $y$, we denote it $Z({\rm sgn},\pi,L,y)$.\end{defn}

It is expected that all the above $L$-functions will eventually 
be shown to be ``essentially
algebraic entire functions'' (Definition 8.5.12 of \cite{go1}). This means
in practice that the $1$-parameter family of power series associated to
the $L$-function, as described above, will have
zeroes in $x^{-1}$ that ``flow continuously;'' this can be
shown for $\zeta_{{\mathcal O}_F}(s)$ (see Theorem 8.9.2 of \cite{go1}).
It is also known for certain $L$-series of Drinfeld modules when
$\A=\Fr[T]$, see \cite{tw1}.

As the zero fields may depend on the sign function, $\pi$ and $y\in \Zp$, it
is important to know how the
zeroes depend on such choices. In fact, we will show in Theorem \ref{indepofall} that
the {\it maximal separable subfields} of the zero fields are independent of
both sgn and $\pi$. Let $L(s)$, $s\in S_\infty$, be a fixed $L$-series of
arithmetic interest as
above. To begin with, let $_i\pi$, $i=1,2$, be two choices of positive
parameters in $\K$ for the {\it same} sign function sgn, and let
$y\in \Zp$. Let $\{_i\lambda^{(y)}_t\}$, $i=1,2$, be the set of zeroes of
$L(x,y)$ as a function in $x$, where the superscript ``$(y)$'' just denotes
the dependence on $y$. One sees from Proposition \ref{dependpi} that to pass
from $\{_1\lambda_{t}^{(y)}\}$ to
$\{_2\lambda_{t}^{(y)}\}$ one just multiplies by $(u^{1/d_{\infty}})^{y}$. Let $_i\K_y$, $i=1,2$, be the extension of
$\K$ generated by $\{_i\lambda^{(y)}_t\}$; that is 
$_i\K_y=Z({\rm sgn},{_i\pi},L,y)$.

\begin{lemma}\label{coeffequal}
Let $\{a_i\}$, $i=1,\ldots n$, be a finite number of elements inside a fixed
algebraic closure $\overline {\mathcal F}$ of a field $\mathcal F$ of 
characteristic $p$. Let $\{u_{i}\}$, $i=1, \ldots, n$, be a collection of
totally inseparable (over $\mathcal F$) elements in $\overline {\mathcal F}$.
Then the
maximal separable subfields of ${\mathcal F}(\{a_i\})$ and 
${\mathcal F}(\{a_iu_{i}\})$ are equal.
\end{lemma}

\begin{proof}
Let $t>0$ be chosen so that $u_{i}^{p^t}\in \mathcal F$ for all $i$. By the multiplicativity
of the inseparability degree over towers of extensions, it is easy to
see that the maximal separable subfield of ${\mathcal F}(\{a_iu_{i}\})$ is the same
as that of ${\mathcal F}(\{a_i^{p^t}\cdot u^{p^t}_{i}\})$. As $u_{i}^{p^t}\in \mathcal F$, this last
field equals ${\mathcal F}(\{a_i^{p^t}\})$. However, the same reasoning as above
shows that the maximal separable subfield of ${\mathcal F}(\{a_i^{p^t}\})$ is the
same as that of ${\mathcal F}(\{a_i\})$ which completes the proof.\end{proof}

The reader will easily see how to generalize Lemma \ref{coeffequal} to
arbitrary collections of elements.

\begin{prop}\label{indofi}
The maximal separable subfields (over $\K$) of $_{i}\K_{y}$, $i=1,2$ are equal.
\end{prop}

\begin{proof}
From Remarks \ref{totins} we see that $u^{1/d_\infty}$ is totally inseparable
over $\K$. The result thus follows from Lemma \ref{coeffequal} and the
remarks made just before Lemma \ref{coeffequal}.
\end{proof}

\begin{rems} \label{opensub} Let $p^{e}$ be the exact power of $p$ dividing
$d_{\infty}$. Let $y\in p^{e}\Zp$. Then Proposition \ref{dependpi} 
immediately shows that the fields $_{i}\K_{y}$, $i=1,2$ are equal. In other
words, there is always a non-trivial open subset of $\Zp$ where the fields
generated by the roots are independent of the choice of positive $\pi$.
\end{rems}

Next we establish that the maximal separable subfield at $y$ of the
root fields is actually independent of the choice of sign function. Let 
${\rm sgn}_{1}$ and ${\rm sgn}_{2}$ be two sign functions. Let
$a\in \A$ and let $_{1}\pi$ be a uniformizer which is positive
for ${\rm sgn}_{1}$ (i.e., ${\rm sgn}_{1}(_{1}\pi)=1$). Thus we can
write 
$$a={\rm sgn}_{1}(a)(_{1}\pi)^{d}u_a\,,$$ where $u_a\in U_1$.
Now let $\zeta={\rm sgn}_{2}(_{1}\pi)$.
Thus
$$a={\rm sgn}_{1}(a)\zeta^{d}\cdot (\zeta^{-1}\cdot {_{1}\pi})^{d}\cdot u_a\,$$
where $_{2}\pi:=\zeta^{-1}\cdot {_{1}\pi}$ is now a positive uniformizer 
for ${\rm sgn}_{2}$. 

\begin{prop}\label{totalind}
Let $\mathfrak I$ be a fractional ideal of $\A$. Then
$\langle {\mathfrak I}\rangle$ is independent of the choice of $_{1}\pi$ or
$_{2}\pi$ as positive uniformizer.
\end{prop}

\begin{proof} Let $\langle {\mathfrak I}\rangle_{i}$, $i=1,2$, be defined with
respect to $_{i}\pi$, $i=1,2$. Now by definition
$\langle (a)\rangle_{i}^{r^{d_{\infty}}-1}=
u_{a}^{r^{d_{\infty}-1}}$ for $i=1,2$. Thus by the unique divisibility of the
$1$-units, we conclude that 
$\langle (a)\rangle_{1}=\langle (a)\rangle_{2}=u_a$. The 
result now follows by the unique divisibility of the $1$-units and the finiteness
of ${\mathcal I}/{\mathcal P}^+$.
\end{proof}

We summarize the above results in the following statement.

\begin{theorem}\label{indepofall}
Let $L(s)$ be an $L$-series of arithmetic interest. Let $y\in \Zp$ be
fixed and let $\K_{y}$ be the extension of $\K$ generated by the zeroes in $x$
of $L(x,y)$. Then the maximal separable subfield of $\K_{y}$ is independent
of the choice of sign function and positive uniformizer.
\end{theorem}

\begin{proof}
This is Propositions \ref{indofi} and  \ref{totalind}. \end{proof}

\begin{defn} \label{zfields} Let $L(s)$ be an $L$-series of arithmetic
interest.\\
1. We denote the common maximal separable subfield of $\K_{y}=Z({\rm sgn},
\pi,L,y)$ by
$Z(L,y)$. This field {\it only} depends on $y\in \Zp$ by Theorem 
\ref{indepofall}. \\
2. We define 
$$Z_{-}(L):=\bigcap_{y\in \Zp}Z(L,y)\,$$
and 
$$Z_{+}(L):=\prod_{y\in \Zp}Z(L,y)\,$$
where product means the compositum. Both of these 
fields are obviously independent
of $y\in\Zp$.\end{defn}

\begin{rems}\label{insep}
Let $\K_1$ be a finite extension of $\K$; so $\K_1$ is also clearly
a local field. By taking roots of a uniformizer, it is easy to see that
any totally inseparable extension of $\K_1$ is uniquely
determined by its degree. If $\K_1$ is assumed to be separable over
$\K$, then one can let $\K_1$ have infinite degree also. In particular,
to find $Z({\rm sgn},\pi,L,y)$ from $Z(L,y)$ one needs only know
the degree of inseparability. As this degree may depend on sgn and
$\pi$ we denote it $i({\rm sgn},\pi,L,y)$. Thus the {\it essential} part
of $Z({\rm sgn},\pi,L,y)$ is {\it precisely} $Z(L,y)$. 
\end{rems}

A complete theory of such $L$-functions
would also give a formula for $i({\rm sgn},\pi,L,y)$. 

Now let $L(s)$ be the $L$-series of a Drinfeld module or general $\A$-module.
Thus we see that $Z(L,y)$ depends {\it only} on the isogeny class
of the Drinfeld module or $\A$-module and $y\in \Zp$.

\begin{questions}\label{ques1} 1. Is it possible to
give a canonical description of $Z(L,y)$ only 
in terms of $y$ and the isogeny class
of the underlying motive?\\
2. Is it possible to give a canonical description of the 
fields $Z_{\pm}(L)$ only
in terms of the isogeny class of the underlying motive?
\end{questions}

We note that, as of this writing,
there seem to be two obvious choices for $Z(L,y)$. The
first is the ``maximal'' choice of the separable closure $\K^{\rm sep}$ of
$\K$ inside ${\mathbf C}_\infty$.
Now the work of Wan, Thakur, Diaz-Vargas, and Sheats, \cite{w1} 
\cite{dv1} \cite{sh1},
shows that the maximal choice is not always correct. The second choice is the
``minimal'' choice which we now define.

To begin with, let $\mathcal F$ be any local non-Archimedean field and 
let ${\mathcal F}_1$ be a finite extension of
${\mathcal F}$. Let $f(x)\in {\mathcal F}_1[[x]]$ be an entire power series with $f(0)=1$. Let
$\{\lambda_t\}$ be the collection of zeroes of $f(x)$ in some fixed
algebraic closure of ${\mathcal F}$ containing ${\mathcal F}_1$. Let ${\mathcal F}_2$ be the extension 
of ${\mathcal F}$ obtained by adjoining to ${\mathcal F}$ the coefficients of $f(x)$, and let
${\mathcal F}_3$ be the extension of ${\mathcal F}$ obtained by adjoining the roots of $f(x)$.

\begin{prop}\label{oneinother}
We have ${\mathcal F}_2\subseteq {\mathcal F}_3$.
\end{prop}

\begin{proof}
Let $t\in \mathbb R$ be such that there exists zeroes of $f(x)$ of
absolute value $t$. Let $f_t(x)=\prod_\lambda (1-x/\lambda)$ where we take
the product of over all $\lambda$ of absolute value $t$ {\it with}
multiplicity. General theory tells us that $f_t(x)$ is a polynomial
with coefficients in ${\mathcal F}_1$. Let ${\mathcal F}_4\subset {\mathcal F}_1$ 
be the extension of ${\mathcal F}$ obtained
by adjoining all the coefficients of all such $f_t(x)$. Note that also
${\mathcal F}_4\subseteq {\mathcal F}_3$ by definition. Now as ${\mathcal F}_1$ is a finite dimensional
extension of ${\mathcal F}$ we also see that ${\mathcal F}_4\subseteq {\mathcal F}_1$ is a finite
dimensional extension of ${\mathcal F}$; thus it is automatically complete. General
theory implies that 
$$f(x)=\prod_t f_t(x)\,,$$
and so the coefficients of $f(x)$ are contained in ${\mathcal F}_4$ by completeness.
As ${\mathcal F}_4\subseteq {\mathcal F}_3$, this completes the proof.\end{proof}

Now let $L(s)$ be our fixed $L$-series as above. By definition, $L(s)$ will
be given as an Euler product over Euler factors of the
form 
$$P_{\mathfrak P}(N{\mathfrak P}^{-s})=P_{\mathfrak P}(x^{-\deg N\mathfrak P}
\langle N{\mathfrak P}\rangle^{-y})\,.$$
Upon expanding out the Euler product, as above, we can write
$$L(s)=\sum_{j=0}^\infty x^{-j}\left( \sum_{\deg N{\mathfrak I}=j}
a_{\mathfrak I}\langle N{\mathfrak I}\rangle^{-y}\right)\,,$$
where $\mathfrak I$ runs over the ideals of the $\A$-integers 
of the base field $F$.

\begin{defn} \label{cfl}
We define $C({\rm sgn},\pi, L,y)$, to be the extension of $\K$ obtained
by adjoining all coefficients 
$$\left\{\sum_{\deg {N\mathfrak I}=j}a_{\mathfrak I}\langle N{\mathfrak I}
\rangle^{-y}\right\}_{j=0}^\infty\,.$$
\end{defn}

\begin{prop} \label{cfieldseq}
The maximal separable (over $\K$) subfield of $C({\rm sgn},\pi,L,y)$ is
independent of the choice of {\rm sgn} and $\pi$.
\end{prop}

\begin{proof}
This follows from Remarks \ref{totins} and Lemma \ref{coeffequal} exactly as in Proposition
\ref{indofi}.
\end{proof}

\begin{rems} \label{containedin}
Note that the proof of Lemma \ref{coeffequal} {\it also} establishes that
the maximal separable subfield of  $C({\rm sgn},\pi, L,y)$ is contained
in the maximal separable subfield of $\K(\{a_{\mathfrak I}\})$. 
Example \ref{cargeo} given below
shows that this containment may be strict.
\end{rems}

\begin{defn}\label{defcl}
1. We denote the common maximal separable subfield of $C({\rm sgn},\pi,L,y)$
by $C(L,y)$. \\2. We set
$$C(L):=\bigcap_y C(L,y)\,.$$
\end{defn}

Since $C(L,y)$ is obviously also independent of any choice of sign function
or uniformizer, we see that this is also true of $C(L)$.

Proposition \ref{oneinother} now immediately implies the next result.

\begin{cor} \label{zinc}
The field $C(L,y)$ is contained in $Z(L,y)$ for all $y\in \Zp$. Thus
$C(L)\subseteq Z_{-}(L)\subseteq Z_{+}(L)$.
\end{cor}

The field $C(L,y)$ is thus the {\it minimal} choice for the fields 
$Z(L,y)$ and $C(L)$ is the minimal choice for $Z_{\pm}(L)$.
It is therefore natural to suspect that $C(L,y)=Z(L,y)$
and $C(L)=Z_{\pm}(L)$. 
However, as in \cite{go1} Section 8.24  one knows that this cannot
always be true and depends on whether the Galois representations that
are used to define the $L$-series are irreducible or not when induced to
${\rm Gal}(k^{\rm sep}/k)$ (if these induced representations remain
irreducible then we call the original motive ``absolutely irreducible''
or ``absolutely simple''). The problem is that
if the representation is not irreducible the $L$-series may factor say into
$L(s)=L_1(s)L_2(s)$ where e.g., $C(L_1,y)$ may be strictly larger than
the original $C(L,y)$ (see the examples  below). Such a factorization
commonly arises from those motives with ``complex multiplication'' by some
order in a
sufficiently large finite extension of $\k$. 

Thus for the moment, let us
assume that our absolutely simple motive has no complex multiplication implying
that $C(L,y)=\K$. Then based on the
evidence available at the moment, as mentioned in the introduction, 
the ``Riemann
hypothesis'' in the context of our $L$-series 
{\it may} turn out to be the statement that
for such irreducible representations $Z(L,y)=C(L,y)=\K$ for all
$y\in \Zp$. Indeed, the results of \cite{w1}, \cite{dv1} and
\cite{sh1} establish that the zeroes of $\zeta_{\Fr[\theta]}(s)$
are both {\it simple} and in $\K$ for {\it any} $y$. However, even
here there is not yet enough
evidence of both a computational and a theoretic nature to know for
certain. For instance, it may turn out that the zeroes are always in
a fixed finite extension of $\K$ (in which case one would like to predict
{\it which} extension this is or at least bound its degree and
discriminant etc.). Future research may perhaps settle this
point. See Remarks \ref{coefffield} for motives with c.m.

\begin{examples}\label{rhirred}
Let $\A=\Fr[T]$, $\k=\Fr(T)$ $k=\Fr(\theta)$ etc. Let $\psi$ be a Drinfeld module over
$k$ with {\it no} complex multiplications over $\overline{k}$. Let $L(s)=
L(\psi,s)$. Then R. Pink
\cite{p1} establishes that the $v$-adic Tate representations of 
${\rm Gal}(k^{\rm sep}/k)$ have {\it open} image. As such they are certainly
irreducible. As the coefficients of the characteristic polynomials
are all in $\A$, one sees immediately that $C(L,y)=\K$ for all
$y$ (so that $C(L)=\K$ also). It is now not obviously 
unreasonable to expect that $C(L,y)=Z(L,y)=\K$ for all $y$.
\end{examples}

When the $L$-series of the motive factors there are still some subtleties
that must be kept in mind. This will be illustrated by the following
examples (still with $\A=\Fr[T]$ etc.) of complex multiplication.
Let $F$ be a finite extension of $k$ with $\psi$ a Drinfeld module
of rank $2$ defined over $F$. We suppose further that $\psi$ has
complex multiplication over $F$ by a separable quadratic extension
$\k_1$ of $\k$ with only {\it one} prime of $\k_1$ lying above $\infty$ --- we also
denote this prime of $\k_1$ by $\infty$. Let $\K_1$ be the
completion of $\k_1$ at $\infty$; thus $[\K_1\colon \K]=2$. Let
$L(\psi,s)$ be the $L$-series of $\psi$ over $F$. As is classically
true, there is a factorization $L(\psi,s)=L(\rho_1,s)L(\rho_2,s)$ into
a product of $L$-series of Hecke characters. Let $L_1(s)$ 
be $L(\rho,s)$ where $\rho$ is either $\rho_1$ or $\rho_2$.
Clearly $C(L,y)=C(L)=\K\subseteq C(L_1,y)\subseteq  Z(L,y)$ for
all $y\in \Zp$. One now needs to
compute $C(L_1,y)$ for such $y$. 
One expects that $C(L_1,y)$ will equal $\K_1$. However,
this is not always true as our first example (Example \ref{carrootuni}) shows.

\begin{examples} \label{carrootuni}
Let $F:={\mathbb F}_{r^2}(\theta)$, $A_1:={\mathbb F}_{r^2}[\theta]$,
$\A_1:={\mathbb F}_{r^2}[T]$ etc. Let
$\psi$ be the Carlitz module for $\A_1$ over $F$ given by
$$T\mapsto \theta \tau^0_1+\tau_1$$
where $\tau_1$ is the $r^2$-power mapping. Of course $\psi$ is the
rank $2$ Drinfeld module for $\A$ with complex multiplication by $\A_1$. The
Hecke character $\rho$ for $\psi$ is the mapping which takes
$(\overline{f})$ where $\overline{f}\in A_1$ is a monic prime to $ f:=
\theta^{-1}(\overline{f})$. Thus 
$$L(\rho,s)=\prod_{\overline{f}}(1-fN(\overline{f})^{-s})^{-1}=
\sum_{\stackrel{\scr \overline{g}\in A_1}{\overline{g}~\rm monic}} 
gN(\overline{g})^{-s}\,.$$
Let $\sigma$ be the non-trivial automorphism of $\A_1/\A$ --- we use the same
notation for the non-trivial automorphism of $A_{1}/A$. Then the above
sum can be written
$$\sum_{\stackrel{\scr \overline{g}\in A}{\overline{g}~\rm monic}}g
N(\overline{g})^{-s}+
\sum_{\stackrel{\scr \{\overline{g},\overline{g}^{\sigma}\}\subset A_1-A}{\overline{g}~\rm monic}}
({ g}+{ g}^\sigma)N(\overline{g})^{-s} \,.$$
Thus in this case $C(L_1,y)=\K$ for all $y$.
\end{examples}

On the other hand, the following example (Example \ref{cargeo}) gives an
example where $C(L_1,y)=\K_1$ for some $y\in \Zp$.

\begin{examples} \label{cargeo}
Let $r=3$, and let 
$\A_1=\A[{\mathbf \lambda}]$ where ${\mathbf \lambda}$ satisfies
${\mathbf \lambda}^2=-T$. 
Let $\psi$ be the Carlitz module over $A[\overline{\lambda}]
=\Fr[\overline{\lambda}]$ given
by $\psi_{\mathbf \lambda}=\overline{\lambda} \tau^0+\tau$ where
$\tau$ is the $r$-th power mapping. Thus $\psi$ is a rank $2$ Drinfeld
module for $\A$ defined over $F=k(\overline{\lambda})=\Fr(\overline{\lambda})$
 given by
$$\psi_T=\theta\tau^0+(\overline{\lambda}^3+\overline{\lambda})\tau+\tau^2=
\theta\tau^0+\overline{\lambda}(1-\theta)\tau+\tau^2\,.$$
Obviously $\psi$ has complex multiplication by $\A_1$. The associated 
Hecke character and $L$-series are just as in Example \ref{carrootuni}.
Note however that the non-trivial automorphism of $\A_1/\A$ does not
take monics to monics. Looking at monics of degree $1$ in $\A_1$ one 
immediately has that for $y\neq 0$ we have $C(L_1,y)=\K_1$. (It is also
very easy to see that $L(\rho,(x,0))\equiv 1$. Thus $C(L_1,0)=\K$ while
there are obviously no zeroes at $y=0$ and, by definition $Z(L_1,0)=\K$.)
 \end{examples}

\begin{rems}\label{coefffield}
The examples just presented
are instances of motives with complex multiplication. We will present here
a variant of the above ``Riemann hypothesis'' that may be more appropriate
in the complex multiplication case.
(See e.g. \cite{ta1} for the  background on the
classical theory --- the characteristic $p$ case
is modeled on this.) Let ${\mathbf F}$ be the complex multiplication field
of an absolutely simple motive $E$ over $\overline{\k}\subset \overline{\K}$; 
thus $\mathbf F$ is realized as
algebraic endomorphisms of the isogeny class of $E$ extending the
canonical $\k$-action.
We will consider the
c.m. field $\mathbf F$ embedded in $\overline{\K}$. Moreover 
{\it once} such an 
embedding is given we also require that if $f\in \mathbf F$ and
$E_{f,\ast}$ is its tangent action on ${\rm Lie}(E)$ then
$$E_{f,\ast}=\overline{f}\cdot{\mathbf 1}+N_f$$
where $\mathbf 1$ is the identity and $N_f$ is nilpotent. This is just
the obvious generalization of the definition of a $T$-module
(e.g., Definition 5.4.5 of \cite{go1}).

Now let $L(s) =L(E,s)$. Note that the coefficients $\{a_{\mathfrak I}\}$
of $L(s)$
are contained in $\mathbf F$ (as we saw in
the above examples). Note also that {\it all the arguments given above
in terms of $\K$ immediately work for $\K_1:={\mathbf F}\cdot \K=\K(
{\mathbf F})$.} Thus,
for instance,
the maximal separable subfield of the extension of $\K_1$ generated by
the zeroes of $L(x,y)$,
for fixed $y$,
is independent of the sign function and positive
uniformizer etc.
Note further that, by Remarks \ref{containedin}, the fields
$C(L,y)$ are all contained in the maximal separable subfield  of $\K_{1}$.
As $\K_1$ obviously does
not depend on $y$,  it {\it may} ultimately be true that for absolutely 
simple motives with c.m., the ``correct'' Riemann hypothesis is the statement
that for fixed $y$ the maximal separable subfield, {\it over
$\K_1$}, of the extension of 
$\K_1$ generated by the zeroes of $L(x,y)$ is precisely $\K_1$. Again, only
future research will determine the truth of this statement. See 
also Question \ref{multicompl} in Section \ref{tmodule}.
\end{rems}

\begin{rems}\label{nofe}
The principle of the 
two $T$'s gives an explanation of why there is no classical style
functional equation for the $L$-functions defined here. Let $\A=\Fr[T]$ etc.
and let $C$ be the Carlitz module with exponential function $e(z)$ where
$z\in C_{\infty}$. Note that, {\it because we are using the two $T$'s 
set-up}, $e(z)\in k[[z]]$. Let $\xi\in C_{\infty}$ be the
period of the Carlitz module so that the lattice of the Carlitz module
is $A\cdot \xi=\A_{\ast}(\xi)$. Let $\Pi(i)\in A$ be the Carlitz factorial. 
The {\it Bernoulli-Carlitz numbers} $BC_{i}$ are defined by
$$\frac{z}{e(z)}=\sum_{i=0}^{\infty}\frac{BC_{i}}{\Pi(i)}z^{i}\,.$$
Thus $\{BC_i\}\subseteq k$.
Because we are using the two $T$'s we see that the definition of
$BC_{i}$ given here is precisely $\theta$ applied to the definition given
on page 354 of \cite{go1}. 

Let $i$ be a positive integer which is divisible by $r-1$.  Note that 
from Definition \ref{zeta} we see that 
$$\zeta_{A}(i):=\zeta_{A}(s_{i})\in \K\,.$$ Thus, from the point of 
view of the two $T$'s, the ``correct'' formulation of Carlitz's 
``Euler'' result is 
$$\zeta_{A}(i)_{\ast}=\xi^{i}\frac{BC_{i}}{\Pi(i)}\,;$$ that is, the 
{\it tangent action at the origin of the multi-valued operator 
$\zeta_{A}(i)$ is precisely scalar multiplication by 
$\xi^{i}\frac{BC_{i}}{\Pi(i)}$.} Equivalently, 
$$\zeta_{A}(i)=\xi^{i}\frac{BC_{i}}{\Pi(i)}\cdot\tau^{0}+ \{{\rm 
higher~terms~in~}\tau\}\,.$$

Now, on the other hand, standard results (8.8.1.1 of \cite{go1}) immediately
imply that $\zeta_{A}(-i)\in \A$. Thus, it is a-priori {\it impossible} to 
match up $\zeta_{A}(-i)$ and $BC_{i}$ etc., in the classical style as
one is an operator and one is just a scalar!
\end{rems}

\begin{questions} \label{higherterms}
Is it possible to give a formula which generalizes the one of Carlitz just
given to the higher terms of the expansion of $\zeta_{A}(i)$ at the
origin?
\end{questions}

Other uses of zeta values at positive integers (see \cite{at1} for instance)
further illustrate that one really uses $\zeta_A(i)_\ast$ as opposed
to the full zeta value. On the other hand, the theory of Drinfeld modules
over finite fields and their zeta functions (\S 4.12 of \cite{go1})
lies squarely in the
universe of $T$ as operator. Indeed, the Frobenius is obviously
an operator on the Drinfeld module and one is interested precisely in the
extension of $\k$ generated by it. Moreover, Gekeler's $Z$-function
(Definition 4.12.24 of \cite{go1}) computes the Euler-Poincar\'e
characteristics $\chi({\mathbb F}_{q^n})$ of the Drinfeld module over the finite fields ${\mathbb F}_{q^n}$, and $\chi({\mathbb F}_{q^n})$ is a principal
ideal in $\A$ etc.

\section{Multi-valued operators and Drinfeld modules}
\label{drinfeld}

Now that we have explained the importance of separability in terms of
zeroes of $L$-series, our goal is to explain its importance in terms of
formal modules and multi-valued operators. In particular we shall eventually 
explain how for general motives it is only possible to extend the formal
module to separable elements; this will be our next section. In this section
we briefly discuss the theory for Drinfeld modules where the result is actually
trivial.

Let $\A$, etc., be general and let $K_1$ be a finite extension of $K$. Let
$\psi$ be a Drinfeld module defined over $K_1$. Thus if
$K_1\{\tau\}$ is the $K_1$-vector space of polynomials in 
$\tau=$ the $r$-th power mapping, then $K_1\{\tau\}$ becomes an $\Fr$-algebra
under composition and $\psi\colon \A\to K_1\{\tau\}$ is
an injection. Let $e(\tau):=e_\psi(\tau)$,
$\log(\tau):=\log_\psi(\tau)$ be the exponential and logarithm of $\psi$
as entire $\Fr$-additive functions over $K_1$.

\begin{defn}\label{multidef1}
A formal power series $f(\tau)$ of the form $f(\tau)=f_\lambda(\tau)
=e(\lambda \log(\tau))$, $\lambda\in C_\infty$,
 is said to be a {\it multi-valued operator} on $\psi$.
\end{defn}

Note that the multi-valued operators clearly form an $\Fr$-subalgebra
of the algebra of all formal power-series in $\tau$.
A multi-valued operator is easily seen to be at most as multi-valued as the
logarithm is; e.g., if $\lambda$ actually is an element of $\A$ then,
of course, $f_{\overline{\lambda}}(\tau)=\psi_\lambda(\tau)$ is an operator on 
$\mathbb{G}_a$. As mentioned in the introduction, one may view a multi-valued
operator as being an operator on {\it all} of $E$ in exactly the same
manner $x^s=e^{s\log(x)}$ is a multi-valued function on 
$\C^\ast$ for any $s\in \C$.

We then have the following simple result.

\begin{prop}\label{drinfeldcase}
The injection $\A \to L\{\tau\}$,
$a \mapsto \psi_a$, extends to an isomorphism of ${\bf C}_\infty$ with
the algebra of multi-valued operators on $\psi$.
\end{prop}

\begin{proof}
By definition, for $a\in \A$ we have
\begin{equation}\label{eqdrinfeldcase} \psi_a(\tau)=e(\overline{a}\log(\tau))\,.
\end{equation}
As there is no obstruction to
replacing $a$ with an arbitrary element in ${\bf C}_\infty$, the result
follows.
\end{proof}

Obviously, the above extension of the $\A$-action to ${\bf C}_{\infty}$ 
is equivalent to extending $\psi$ to a formal ${\bf C}_{\infty}$-module.
\section{Separability, Multi-valued operators and general $T$-modules}\label{tmodule}

We now discuss the applications of separability to arbitrary motives. Due to
the lack of a worked out theory of uniformizability of motives for general
$\A$, we restrict ourselves to $\A=\Fr[T]$ etc. However, it is clear that once this
theory is worked out, the arguments given here will immediately apply in
general.
We begin this section by showing that Proposition \ref{drinfeldcase} is
profoundly misleading in general.

\begin{examples}\label{exq2}
Let $E=C^{\otimes 2}$ be the second tensor power of the Carlitz module. For each $a\in \A$ we let $C_a^{\otimes 2}$ be the
associated $\Fr$-linear algebraic endomorphism of ${\mathbb G}_a^2$ and we let
$C_{a,\ast}^{\otimes 2}$ be the associated action on the tangent space at the
origin. Thus, for $a=T$ one has
\begin{equation}\label{Tstar}
C^{\otimes 2}_{T,\ast}=\left(\begin{array}{cc}\theta& 1\\ 0& \theta
\\\end{array} \right)\,.\end{equation}
Now if it were possible in general to extend $C^{\otimes 2}_a$ to
{\it all} $a\in {\bf C}_\infty$ then
$C^{\otimes 2}_{T^{1/2},\ast}$ would be a solution to the matrix equation
$x^2=C^{\otimes 2}_{T,\ast}$ where $C^{\otimes 2}_{T,\ast}$ is given as in
Equation \ref{Tstar}. However, it is easy to see that this equation is
inconsistent for $p=2$.
\end{examples}

The reader will easily see how to construct other such examples. Indeed,
let $N$ be a nilpotent matrix. Then, of course, if $x=x\cdot {\bf 1}$ is the
scalar matrix, one has
$$(x+N)^{p^t}=x^{p^t}+N^{p^t}\,.$$
Thus for $t\gg 0$, $(x+N)^{p^t}$ is also a scalar matrix and from this
other examples are readily found.

Now let $E$ be an arbitrary uniformizable $T$-module defined over
$C_\infty$ and where the 
underlying algebraic group is isomorphic to $\mathbb{G}_a^t$ (i.e.,
$E$ has ``dimension $t$''). We note that
Anderson's theory implies that {\it any} tensor product of Drinfeld modules
is automatically uniformizable. Let $e({\bf z}):=e_E({\bf z})$ and
$\log({\bf z}):=\log_E({\bf z})$ be the exponential and logarithm of
$E$ with respect to some coordinate system on the underlying
algebraic group; thus $\bf z$ is a vector (and so that uniformizability is
exactly equivalent to $e({\bf z})$ being surjective). 
Let $\tau=\tau({\bf z})$ be the mapping that raises each coefficient
of $\bf z$ to the $r$-th power. As the exponential and logarithm
functions are $\Fr$-linear, they can 
be written as power series in $\tau$ which we denote by
$e(\tau)$, $\log(\tau)$ etc. More precisely each of these two
functions may be written as $\sum_{i=0}^\infty c_i\tau^i$
were the $c_i$ are $t\times t$ matrices with coefficients in the smallest
extension of $k\subset C_\infty$ which contains the coefficients of the
$T$-action.

\begin{defn}\label{multidef2}
Let $f(\tau)$ be any power
series in $\tau$ of the form $\sum_{i=0}^\infty a_i\tau^i$ where
the $a_i$ are $t\times t$ matrices with $C_\infty$-coefficients. We say
that $f({\tau})$ is a {\it multi-valued operator on $E$} if and only if it
may be written in the form $f(\tau)=e(M\log(\tau))$ where $M$ is a 
$t\times t$ matrix with $C_\infty$-coefficients.
\end{defn}

Definition \ref{multidef2} is the obvious generalization of 
Definition \ref{multidef1}. Moreover, it is clear that such 
multi-valued operators form a sub-algebra of the algebra of
all formal power series in $\tau$ and that this algebra is
isomorphic to the $\Fr$-algebra of $t\times t$ matrices over
$C_\infty$. 

It is our goal in this subsection to show that, in spite of
Example \ref{exq2}, the $T$-action on $E$ still
has a vast extension to the realm of multi-valued operators. More 
precisely we will show the following result.

\begin{theorem}\label{mainth}
The $T$-action of $E$ can be uniquely extended to an injection of
$\Fr$-algebras,
$\lambda \mapsto {E}_{\lambda}$, of $\K^{\rm sep}$ into the algebra of
multi-valued operators. Moreover, if ${E}_{{\lambda},\ast}$
is the induced action on the tangent space at the origin, then
$${E}_{{\lambda},\ast}=\overline{\lambda}\cdot{\bf 1}+
N_{\lambda}=\overline{\lambda} +N_{\lambda}$$ 
where
$N_\lambda$ is nilpotent.
\end{theorem}

As in Proposition \ref{drinfeldcase}, through the use of the
exponential and logarithm of $E$ it is clear that
the result will be proved if we can establish that the tangent action
$${E}_{T,\ast}=\theta\cdot {\bf 1}+N_T=\theta+N_{T}$$
of $E_T$
uniquely extends to an injection ${E}_{x,\ast}$
of $\K^{\sep}$ into the $\Fr$-algebra of $t\times t$ matrices with
coefficients in $C_\infty$ with the required form. That is, we may work
``infinitesimally'' with nilpotent matrices.

Theorem \ref{mainth} will follow immediately from our next two
propositions. Let us set $\epsilon:=N_T$ where $N_T$ is the nilpotent
in $E_{T,\ast}$.

\begin{prop}\label{toK} The $T$ action on $E$ extends uniquely to
an injection of $\K$ into the algebra of multi-valued operators
such that if $x\in \K$, then $E_{x,\ast}=\overline{x}+N_x$ 
and $N_x$ is a polynomial in $\epsilon$.
\end{prop}

\begin{proof} The last part is clear for all $a\in \A$. Now let us
examine what happens for $a=1/T$. As $E_{T,\ast}=\theta\cdot{\bf 1} + \epsilon
=\theta +\epsilon$, we are
led to define
$$E_{1/T,\ast}:=(\theta+\epsilon)^{-1}=\frac{1}{\theta}\cdot\frac{1}{1+\epsilon/\theta}=
\frac{1}{\theta}\cdot(1-\epsilon/\theta+\epsilon^2/\theta^2+\cdots...)\,.$$
As $\epsilon$ is nilpotent, so is $\epsilon/\theta$ and the result follows for $a=1/T$. Continuity
now finishes the proof for all $x\in \K$.
\end{proof}

Suppose now that $\epsilon$ satisfies the minimal polynomial $u^t=0$ and
let $_1C_\infty:=C_\infty[\epsilon]$ be the $t$-th order infinitesimal thickening
of $C_\infty$. Let $x\in \K$. 
By the proposition, we can represent 
 $E_{x,\ast}$ as an element of $_1C_\infty$; we will denote
the nilpotent part of $E_{{x},\ast}$ by $\epsilon_{x}$.
Let $0\neq {\lambda}\in \K^{\sep}$ and let $\lambda$ satisfy
the separable polynomial equation $f(u)=\sum_{i=0}^e{a}_iu^i=0$ with
$\K$-coefficients  and ${a}_0\neq 0$. 
If we can find an extension of the injection
${x}\mapsto E_{{x},\ast}$ to $\lambda$, then 
clearly $E_{{\lambda},\ast}=
\overline{\lambda}+\epsilon_{\lambda}=\overline{\lambda}+\epsilon_{\lambda}$
satisfies the equation {\it with coefficients in $_1C_\infty$}
\begin{equation}\label{basiceq}
f_1(u):=\sum_{i=0}^e E_{{a}_i,\ast}u^i=\sum_{i=0}^e
(\overline{a_i}+\epsilon_{{a}_i})u^i=0\,.
\end{equation}
Let $\dis f_2(u):=\sum_{i=1}^e(\overline{a_i}+\epsilon_{a_i})u^i+
\overline{a_0}$\,; so that now $f_2(\overline{\lambda}+
\epsilon_\lambda)$ equals $-\epsilon_{{a}_0}$ and is nilpotent.
Set $f_3(u):=f_2(\overline{\lambda}+u)$; thus we deduce that
\begin{equation}\label{2basiceq}
f_3(\epsilon_{\lambda})=-\epsilon_{{a}_0}\,,\end{equation}
and that $f_3$ takes the origin to the origin in ${\rm Spec}(_1C_\infty[u])$.
Note that $f_3^\prime(0)\equiv \overline{f^\prime({\lambda})}\pmod{
(\epsilon)}$. As $f^\prime({\lambda})\neq 0$, $f_3(u)$ is  also \'etale at
the origin. 

\begin{prop}\label{toKbar}
There exits a unique solution $\epsilon_{\lambda}\in {_1C_\infty}$ to {\rm Equation
\ref{2basiceq}}.\end{prop}

\begin{proof}
Since the map $f_3(u)$ is \'etale, there is a formal inverse to it. 
As we are dealing with nilpotents, we may 
evaluate it at $-\epsilon_{{a}_0}$ to find $\epsilon_\lambda$ and to see that
it lies in $C_\infty[\epsilon]$. 
\end{proof}

Theorem \ref{mainth} is now easily seen. Indeed, Proposition \ref{toKbar} certainly
gives the extension to $\K({\lambda})$ and the uniqueness implies that
we can actually piece together these liftings all the way up to $\K^{\sep}$.

\begin{rems} The argument in Proposition \ref{toKbar} shows that the
only possible lifting of $E_{{x},\ast}$ to the algebraic closure of
$\Fr$ must be as scalars; i.e., the nilpotent part vanishes.
\end{rems}

\begin{examples}\label{epsilonlambda}
We put ourselves back in the situation of $E=C^{\otimes 2}$ as
in Example \ref{exq2}. Set $r=5$. We will show how our method can be used to lift
the definition of $C^{\otimes 2}_{x}$
to $\lambda$ satisfying the $\K$-equation in $u$
\begin{equation}\label{lambdaeq} u^{2}+Tu+T^{3}=0\,.\end{equation}
It is easy to see that this equation is irreducible and separable over $\K$.
We will use the notation just given above; so $C^{\otimes 2}_{T,\ast}=\theta
+\epsilon$ and $\epsilon^{2}=0$. Thus $C^{\otimes 2}_{{\lambda},\ast}=
\overline{\lambda}+
\epsilon_{\lambda}$ where $\epsilon_{\lambda}$ needs to be found. 
Using  Equation \ref{lambdaeq} we see that
$$(\overline{\lambda}+\epsilon_{\lambda})^{2}+(\theta+\epsilon)
(\overline{\lambda}+\epsilon_{\lambda})+
\theta^{3}+3\theta^{2}\epsilon=0\,,$$
or
$$\epsilon_{\lambda}^{2}+(2\overline{\lambda} +\theta +\epsilon)
5\epsilon_{\lambda}+
({\overline{\lambda}}^{2}+\theta\overline{\lambda}+\theta^{3})+(
\overline{\lambda}+3\theta^{2})\epsilon=0\,.$$
As ${\overline{\lambda}}^{2}+\theta\overline{\lambda}+\theta^{3}=0$, we see that
$$\epsilon_{\lambda}^{2}+(2\overline{\lambda}
+\theta+\epsilon)\epsilon_{\lambda}=
-(\overline{\lambda}+3\theta^{2})\epsilon\,.$$
As $\epsilon^{2}=0$, one easily solves to see that
$$\epsilon_{\lambda}=\frac{-(\overline{\lambda}+3\theta^{2})\epsilon}
{2\overline{\lambda}+\theta+
\epsilon}=\left(\frac{-\overline{\lambda}-3\theta^{2}}{2
\overline{\lambda}+\theta}\right)\epsilon\,.$$
\end{examples}

We now give another approach to Theorem \ref{mainth} through the use of
differential calculus in characteristic $p$. To see how this can be accomplished
we continue to examine the second tensor power of the Carlitz module.

\begin{examples}\label{exgenq}
(Example \ref{epsilonlambda} redux.)
It is easy to see that for ${a}\in \A$ one has
\begin{equation}\label{astar}
E_{{a},\ast}=\left(\begin{array}{cc}\overline{a}& \overline{a^\prime}\\ 0& 
\overline{a}
\\\end{array} \right)\,,\end{equation}
where $\overline{a}=\theta(a)$, $\dis \overline{a^{\prime}}=
\overline{\frac{d{a}}{dT}}=\frac{d\overline{a}}{d\theta}$ etc. 
It is elementary to see that
injection $a\mapsto E_a$ extends to $a\in \K$ in {\it exactly} the same
form. Now, as is very well-known,
the derivation $x\mapsto \frac{dx}{dT}$ on  $\A$
has a unique extension to
$\K^{\sep}$. Thus this extended derivation in turn furnishes the desired extension
of $C^{\otimes 2}_{{x},\ast}$ to all separable elements and, in particular,
$\lambda$ of Example \ref{epsilonlambda}. In fact, to find 
$\theta(\frac{d{\lambda}}{dT}):=\frac{d\overline{\lambda}}{d\theta}$
one first differentiates implicitly Equation \ref{lambdaeq} to find
$$(2{\lambda}+T)\frac{d{\lambda}}{dT}+({\lambda}+3T^{2})=0\,$$
and then applies $\theta$ after solving for
$\frac{d{\lambda}}{dT}$. The answer is now easily seen to agree with that of
Example \ref{epsilonlambda}.
\end{examples}

The above example can be extended greatly. Let $n$ and $j$ be nonnegative
integers and let $\binom{n}{j}$ be the usual binomial coefficient. Let
${\mathfrak D}_{j}$ be the differential operator on $\A$ which has
$${\mathfrak D}_{j}T^{n}=\binom{n}{j}T^{n-j}\,.$$
(N.B.: in characteristic $0$,  and {\it only} characteristic
$0$, one has ${\mathfrak D}_{j}=
\frac{{\mathfrak
D}_{1}^{j}}{j!}$.) The operators $\{{\mathfrak D}_{j}\}$ are called ``hyperderivatives''
(or ``higher derivatives'') and form the basis of characteristic $p$ calculus.
They arose first in the paper of Hasse and Schmidt \cite{hs1}.
It is simple to see that the collection $\{{\mathfrak D}_{j}\}_{j=0}^{\infty}$
satisfies the ``Leibniz identity:''
$${\mathfrak D}_{n}(uv)=\sum_{i=0}^{n}{\mathfrak D}_{i}(u){\mathfrak D}_{n-i}(v)\,.$$
If now $\epsilon$ is a nilpotent with $\epsilon^{n}=0$, then, as is well known,
the Leibniz formula implies that the map from $\A$ to $\A[\epsilon]$
given by 
$$a\mapsto \sum_{i=0}^{n}{\mathfrak D}_{i}(a)\epsilon^{i}$$
is a homomorphism.

Let $E$ now be an arbitrary $T$-module of dimension $t$ over a finite extension $K_{1}$ of
$K$. As above let $\epsilon:=N_{T}$. 
As the elements of $K_{1}$ are {\it scalars}, the next result follows
immediately.

\begin{prop}\label{conhyper}
Let $a\in \A$. Then $$E_{a,\ast}=\overline{\left(\sum_{i=0}^{t}{\mathfrak D}_{i}(a)
\epsilon^{i}\right)}\in A[\epsilon]\,.$$
\end{prop}

Thus we see that Theorem \ref{mainth} immediately implies the liftability
of hyperderivatives to arbitrary separable extensions of $\K$; in fact
the techniques of Proposition \ref{toKbar} can be also used to show this
well-known fact (see, e.g., \cite{ikn1}) in complete generality. Conversely, Theorem \ref{mainth}
follows directly once we know the liftability of such operators.

The results of Sections \ref{lfunction} and \ref{tmodule} now imply the following
basic question. Let $L(M,s)$ be an $L$-function of arithmetic interest associated
to a motive $M$ defined over a finite extension $L$ of $k$.
Recall that, as in Remarks \ref{zetaL}, $L(M\otimes
C^{\otimes n},s)=L(M,s-n)$. Let $\rho:L\to
\overline{K}$ be an embedding. Thus via $\rho$ we obtain a motive $M^{\rho}$ defined
over a finite extension of $K$ which we may assume is
uniformizable (e.g., if $M$ is any tensor product of Drinfeld modules).
We can therefore speak about multi-valued operators on $M^{\rho}$.

\begin{questions}\label{bigbasic}
Is it possible to canonically describe the fields $C(L,y)$, $Z(L,y)$ etc.
as subalgebras of the algebra of multi-valued 
operators on $M^{\rho}$ and $M^{\rho}\otimes C^{\otimes n}$?
\end{questions}

For instance, let $\psi$ be a Drinfeld module defined over 
$k={\mathbb F}_r(\theta)$. Obviously, there is only one embedding
into $\overline{K}$. One can then look for descriptions of the 
maximal separable (over $\K$) subfields of the zero fields
of $L(\psi,s)$ in terms of multi-valued operators on $\psi$ {\it and}
the motives $\psi\otimes C^{\otimes n}$ where $n\geq 1$.

As in Remarks \ref{coefffield} we can refine the above statements when
the motive has complex multiplication. Before giving this refinement there
is a subtlety that must first be discussed. Let $E$ be
a motive with complex multiplication 
by a field $\mathbf F$ considered as a subfield of
$\overline{\K}$. As in Remarks \ref{coefffield}, we require that
if $f\in \mathbf F$ then then tangent action, $E_{f,\ast}$, is of
the familiar form $E_{f,\ast}=\overline{f}+N_f$ where $N_f$ is nilpotent.
The subtlety is that we will need the same statement to be true of the
canonical action of $\mathbf F$ on all tensor powers 
$E\otimes C^{\otimes n}$. Our next example will show the reader just
how easily this may be established.

\begin{examples}\label{cargeoredux}
Let $r=3$ and $\A_1$, $\psi$ etc., be as in Example \ref{cargeo}.
We set ${\mathbf F}=\k(\lambda)$ where $\lambda^2=-T$; so $\mathbf F$ is
just the quotient field of $\A_1$. We consider $\mathbf F$ as
embedded in $\overline{\K}$. The $T$-module $\psi\otimes C$ is
constructed as in Section 5.7 of \cite{go1}. That is, we first pass to the
perfection $F_1:=F^{\rm perf}$
of the field $F=k(\overline{\lambda})$. Next we let
$M_\psi$ and $M_C$ be the $T$-modules associated to $\psi$ and $C$.
By definition both of these are modules over the ring
$F_1[T]$. Then $\psi\otimes C$ is the $T$-module associated to the
module $M_1:=M_\psi\otimes_{F_1[T]}M_C$ equipped with the {\it diagonal}
action of $\tau$ (= the $r$-th power mapping). That is 
$$\tau(a\otimes b):=\tau(a)\otimes \tau(b)\,.$$
The tangent space of $\psi\otimes C$ is dual to $M_1/\tau M_1$ (Lemma
5.4.7 of \cite{go1}), and so we are reduced to examining the action 
of an element $f\in {\mathbf F}$ on the isogeny class
of $M_1/\tau M_1$. Now without loss of
generality we can assume that $f$ gives rise to 
an actual action on $M_1$ itself.
Note that $T-\theta$ kills both $M_\psi/\tau M_\psi$ and $M_C/\tau M_C$
(as these arise from Drinfeld modules which are $1$-dimensional), and so
$$(T-\theta)^2M_1/\tau M_1=\{0\}\,.$$
Now $\lambda^2+T=0$, which implies
$$(x-\lambda)(x+\lambda)=x^2+T\,.$$
Thus 
$$(\overline{\lambda}-\lambda)(\overline{\lambda}+\lambda)=T-\theta\,.$$
Consequently, by what we just saw,
$(\overline{\lambda}-\lambda)^2(\overline{\lambda}+\lambda)^2$
kills $M_1/\tau M_1$. But the tangent action of $\lambda +\overline{\lambda}$
is clearly invertible and this implies
that $(\lambda-\overline{\lambda})^2$ kills
$M_1/\tau M_1$ as required.
\end{examples}

We can now turn to the version of Question \ref{bigbasic} appropriate
for the case of complex multiplication.
Suppose we have a motive $M$ with complex multiplication as in
Remarks \ref{coefffield}. Let $\mathbf F$ be the complex multiplication
field and $\K_1={\mathbf F}\cdot \K$ as before. As the tangent
action of elements of $\mathbf F$ is assumed to have the same form
as the tangent 
action of $T$ (i.e., ``scalar multiplication + nilpotent''), it is
immediate that the tangent action continues to elements of $\K_1$. Then
the above arguments immediately show that it also extends to 
the separable closure $\K_1^{\rm sep}$ of $\K_1$ inside $\overline{\K}$.
Moreover we have also seen
that $L(M,s)$ has coefficients in $\K_1$.

\begin{questions}\label{multicompl}
Is it possible to canonically describe the maximal separable 
(over $\K_1$) subfields of the zero
fields of $L(M,s)$ as multi-valued operators on $M$ and $M\otimes C^{\otimes n}$ 
etc.?
\end{questions}

\section{The $v$-adic and global theory} \label{vadic}
In this last section we will briefly discuss the $v$-adic and global
versions of
the above results. Let $v\in {\rm Spec}(\A)$ be fixed and let $\k_{v}$ be the
completion of $\k$ at $v$. Let $\overline{\k_{v}}$ be a fixed algebraic
closure of $\k_{v}$ equipped with its canonical metric. Let ${\mathbf C}_{v}$ be
the completion of $\overline{\k_{v}}$.

As in Section \ref{lfunction}, for an ideal $\mathfrak I$ of $\A$ one
has the element ${\mathfrak I}^{s_{1}}\in \overline{\k}\subset
{\mathbf C}_\infty$. The collection of
all such elements generates an extension $\mathbf V$ of $\k$ which can
easily be shown to be finite (Prop. 8.2.9 of \cite{go1}). Let $\sigma \colon
{\mathbf V}\to \overline{\k_{v}}$ be an embedding over $\k$ and let
$\k_{\sigma,v}$ be the compositum of $\sigma({\mathbf V})$ and $\k_{v}$. This
extension is finite over $\k_{v}$ and has $\A_{v}$-integers $\A_{\sigma,v}$.
Let $d_{v}$ be the degree of $v$ over $\Fr$, and let $f_{\sigma}$ be the
residue degree of $\A_{\sigma,v}$ over $\A_{v}$.

\begin{defn}\label{1vadicexp}
 We set 
$$S_{\sigma,v}:=\lim_{\stackrel{\longleftarrow}{j}} \Z/
(p^{j}(p^{d_{v}f_{\sigma}}-1))\isom \Zp\times \Z/(r^{d_{v}f_{\sigma}}-1)\,.$$ 
\end{defn}

Note that $S_{\sigma,v}$ is obviously a ring.

Let $\beta\in\A_{\sigma,v}^{\ast}$. Then, as is completely standard,
$\beta$ has a canonical decomposition $$\beta=\omega_{\sigma,v}(\beta)
\langle \beta
\rangle_{\sigma,v}$$ where $\langle \beta \rangle_{\sigma,v}$ is a $1$-unit
and $\omega_{\sigma,v}$ is the Teichm\"uller representative. Let
$s_{v}=(s_{v,0},s_{v,1})\in S_{\sigma,v}$, where
$s_v$ is not to be confused with $s_j\in S_\infty$. We then set
$$\beta^{s_{v}}:=\omega_{\sigma,v}(\beta)^{s_{v,1}}
\langle \beta\rangle_{\sigma,v}^{s_{v,0}}\,.$$

Now let $\mathfrak I$ be any ideal of $\A$ which is prime to $v$. 
By definition,
if $o$ is the order of ${\mathcal I}/{\mathcal P}^{+}$
 then ${\mathfrak I}^{o}=(i)$ 
where $i$ is positive and prime to $v$. Thus $({\mathfrak I}^{s_{1}})^{o}=i$
which implies that ${\mathfrak I}^{s_{1}}$ is also prime to $v$. 
Our next definition is
then the $v$-adic version of Definition \ref{exponen}.

\begin{defn}\label{vexponen}
Let $\mathfrak I$ be as above and let $e_{v}:=(x_{v},s_{v})\in 
{\mathbf C}_{v}^{\ast}\times S_{\sigma,v}$. Then we set
$${\mathfrak I}^{e_{v}}:={\mathfrak I}^{(x_{v},s_{v})}:=
x_{v}^{\deg \mathfrak I}\sigma({\mathfrak I}^{s_{1}})^{s_{v}}\,.$$
\end{defn}

It is now a very simple matter to use Definition \ref{vexponen} to
define the {\it $v$-adic $L$-functions} of a motive $M$ as an Euler 
product over primes {\it not} dividing $v$ --- this $L$-function is denoted
$L_{\sigma,v}(M,e_{v})$ etc. For the details we refer the reader to
\cite{go1}, Definition 8.6.8.
One can also define the zero fields  of $L_{\sigma,v}(M,e_{v})$ exactly as in Definition
\ref{zerofield}. This field will depend on the original motive $M$,
$v\in {\rm Spec}(\A)$, the sign function sgn, the root of unity $\zeta$ of
Proposition \ref{idependpi}, and $s_{v}\in S_{v}$. 

\begin{examples} \label{vadiczeta}
Let $\A=\Fr[T]$, $A=\Fr[\theta]$ and
$$\zeta_{A}(s)=\sum_{n~\rm monic}n^{-s}=
\prod_{f~\rm monic~ prime}(1-f^{-s})^{-1}\,,$$
as in Example \ref{zeta}. Note that as $\A$ has class number $1$, one
immediately deduces that ${\mathbf V}=\k$; we will thus drop any reference
to the obviously canonical embedding $\sigma\colon \k \to \k_{v}$ for fixed
$v\in {\rm Spec}(\A)$. 
In terms of Euler products, by definition we have
$$\zeta_{A,v}(e_{v})=\prod_{\stackrel{\scr f~\rm monic~prime}{\scr 
f\not \equiv 0~{\rm mod}\, v}}(1-f^{-e_{v}})^{-1}\,$$
etc. This Euler product converges for those $e_{v}$ with $|x_{v}|_{v}> 1$.
Upon expanding out the product, we have
$$\zeta_{A,v}(e_{v})=\sum_{\stackrel{\scr n~{\rm monic}}{n\not \equiv 0~ {\rm mod}\,{v}}}
n^{-e_{v}}=\sum_{\stackrel{\scr n~{\rm monic}}{n\not \equiv 0~ {\rm mod}\,{v}}}
x_{v}^{-\deg n}n^{-s_{v}}\,.$$
By grouping together according to degree, we see
$$\zeta_{A,v}(e_{v})=\sum_{j=0}^{\infty}x_{v}^{-j}\left(
\sum_{\stackrel{\scr n~{\rm monic}}{\stackrel{\scr \deg n=j}{\scr n \not \equiv
0~{\rm mod}\, v}}}n^{-s_{v}}\right)\,.$$
Basic estimates (Lemma 8.8.1 of \cite{go1}) imply that
for fixed $s_{v}\in S_{v}$, $\zeta_{A,v}(x_{v},s_{v})$ converges to an
entire $v$-adic power series in $x_{v}^{-1}$ such that the zeroes ``flow
continuously'' in $s_{v}$.
\end{examples}

There is another way to obtain these $v$-adic functions via
``essential algebraicity.'' For our purposes
here, we will illustrate this by continuing to discuss
$\zeta_{A,v}(e_v)$ where $A=\Fr[\theta]$; in general, we refer the
reader to Sections 8.5 and 8.12 of \cite{go1}.

\begin{examples}\label{vadiczetacont}
We continue here with the set-up of Example \ref{vadiczeta}. Let $\pi$
be our fixed positive uniformizer.
We set $z_\zeta(x,0)\equiv 1$, and for $j$ a positive integer,
$$z_\zeta(x,-j):=\zeta_A(x\pi^j,-j)=\sum_{i=0}^\infty x^{-i}
\left(\sum_{\stackrel{\scr n~\rm monic}{\scr \deg n=i}}n^j\right)\,.$$
As $\zeta_A(s)$ is an entire function on $S_\infty$, the power series
$z_\zeta(x,-j)$ is entire in $x^{-1}$. But as the coefficients are in $\A$, one
sees immediately that this forces $z_\zeta(x,-j)$ to be a {\it polynomial}
in $x^{-1}$. Moreover, if $j$ is positive and divisible by $r-1$, then there is
the {\it trivial zero}
$$\zeta_A(-j)=\zeta_A(s_{-j})=z_\zeta(1,-j)=0\,.$$
Thus we set $\tilde{z}_\zeta(x,-j)=z_\zeta(x,-j)$ if $j\not \equiv
0~{\rm mod}\,(r-1)$, and $\tilde{z}_\zeta(x,-j)=z_\zeta(x,-j)/(1-x^{-1})$ 
for those positive $j\equiv 0~{\rm mod}\,(r-1)$.

The polynomials $\{z_\zeta(x,-j)\}$ and $\{\tilde{z}_\zeta(x,-j)\}$,
are obviously {\it global objects} as they have coefficients in $\A$.
Their importance is precisely the following. Let $v\in {\rm Spec}(\A)$
and substitute $x_v$ for $x$. Then the polynomials
{\it $\{z_\zeta(x_v,-j)\}$ interpolate $v$-adically to the entire functions
$\zeta_{A,v}(e_v)$.} Thus, much information about these $v$-adic functions
can be obtained from $\zeta_\zeta (x,-j)$. Note also that
the process of interpolation removes the Euler factor at $v$ from
$z_\zeta(x_v,-j)$.
\end{examples}

It is certainly very reasonable to inquire about a possible $v$-adic
version of the Riemann hypothesis of \cite{w1}, \cite{dv1} and \cite{sh1}. Our
first result along these lines is contained in a clever observation of
D. Wan that {\it some} $v$-adic information can actually be gleaned from
results already established at $\infty$. This is contained in the
next result.

\begin{prop} \label{wans}
Let $v$ be prime of degree $1$ in $\A=\Fr[T]$. Let $j$ be an element of the 
ideal $(r-1)S_v$. Then the zeroes (in $x_v$) of $\zeta_{A,v}(x_v, j)$ are
simple and in $\k_v$.
\end{prop}

\begin{proof} The idea of the proof is to exploit the isomorphism between
$\K=\k_\infty$ and $\k_v$ (remember $\deg v=1$!). Without loss of generality,
we can set $v=(T)$. Furthermore, we begin by letting $j$ be a positive
integer divisible by $r-1$ and we {\it choose} our positive uniformizer
to be $\pi=1/T$.
Now the coefficient of $x^{-d}$ in $\zeta_A(x,-j)$ is precisely the
sum of $\langle n\rangle^j$ where $\deg n=j$ and $n$ is monic, while the
coefficient of $x_v^{-d}$ in $\zeta_{A,v}(x_v,-j)$ is the sum of
$n^j$ such that $n$ is monic of degree $d$ and $n\not \equiv 0~{\rm mod}\,v$.
This last condition is the same as saying that $n$ has non-vanishing
constant term. 

The set $\{\langle n \rangle\}$, where $n$ is monic, ranges over
all polynomials $f(1/T)$ in $1/T$ with constant term $1$ and degree (in $1/T$) 
$<d$. Moreover, as $j$ is divisible by $r-1$, the set
$\{\langle n \rangle^j\}$ is the same as the set
$f(1/T)^j$ where $f(u)$ is a monic 
polynomial of degree $<d$ and has {\it non-vanishing}
constant term. 

Let us denote by $\zeta_{A,v}(x,-j)$ the function obtained by replacing
$x_v$ by $x$ in $\zeta_{A,v}(x_v,-j)$ and applying the isomorphism
$\k_v \to \k_\infty$ given by $T\mapsto 1/T$. The above now implies
that
$$(1-x^{-1})^{-1}\zeta_{A,v}(x,-j)=\zeta_A(x,-j)\,.$$
The result for positive $j$ divisible by $r-1$ follows immediately when
we recall that $\zeta_{A,v}(e_v)$ has an Euler product over all
primes {\it not} equal to $v$. 
The general result then follows simply by passing to the limit.
\end{proof}

\begin{rems} \label{roberts1}
Proposition \ref{wans} seems to indicate that there may also be a Riemann
hypothesis type phenomenon for the finite primes of $\A$. But the situation
is definitely more subtle than the proposition indicates. Indeed, the
work of Wan, Diaz-Vargas and Sheats establishes that the polynomials
$z_{\zeta}(x,-j)$ of Example \ref{vadiczetacont} are {\it separable}. Now
if Proposition \ref{wans} was true for all $j$ and {\it all} finite $v$,
then these polynomials would have to split totally in $\k$. However,
recent computer work by J. Roberts seems to indicate that the polynomials
$\tilde{z}_{\zeta}(x,-j)$ are {\it irreducible}. Indeed, in the examples
computed, Eisenstein's irreducibility criterion works. Therefore, if there
is a $v$-adic Riemann hypothesis in general, one should expect non-trivial
finite (at least) extensions to be involved. As always, a complete theory
would predict these extensions etc.
\end{rems}

\begin{questions} \label{roberts2}
Are the polynomials $\tilde{z}_{\zeta}(x,-j)$ irreducible in general?
\end{questions}

\begin{rems} \label{roberts3} 
Using standard results in algebra, the Galois groups of some
of the polynomials $\tilde{z}_{\zeta}(x,-j)$ have been computed where
$$\deg_{x^{-1}}(\tilde{z}_{\zeta}(x,-j))\leq 4\,.$$
These
have all turned out to be the {\it full symmetric group} and so, in particular,
are non-abelian in general. Even for small $r$ and $j$, the polynomials in
$T$ that appear as coefficients are quite large and a great deal of computer
time was necessary. For instance, for $r=5$ and $j=1249$, the degree in
$x^{-1}$ of $z_{\zeta}(x,-j)$ is $4$. The constant term is $1$, of course,
the coefficient of $x^{-1}$ has degree $1245$ in $T$, the coefficient
of $x^{-2}$ has degree $2470$, the coefficient of $x^{-3}$ has degree
$3595$, and the coefficient of $x^{-4}$ has degree $4220$. The resolvent
cubic was shown to be irreducible modulo the sixth degree prime
$T^{6}+T^{5}+T^{4}+T^{3}+T^{2}+T+1$ by brute force. The discriminant
of this cubic was computed and seen not to be a square thereby giving the
result. 
\end{rems}

\begin{questions} \label{roberts4}
Is the Galois group of $\tilde{z}_{\zeta}(x,-j)$ the full symmetric group
in general?
\end{questions}

Let $\k(j)\subset {\mathbf C}_\infty$ be the splitting field of
$\tilde{z}_\zeta(x,-j)$. As $\deg \infty =1$, it is easy to see that
this field does not depend upon the choice of sign function and only depends
on $j$, etc.
One would like to also know what are the primes, besides
$\infty$ (and perhaps those of degree $1$), that split completely in 
$\k(j)$ as well as the discriminant (to $\k$) etc. 

Finally, let us define $\k_\zeta$ to be the compositum of all such
$\k(j)$. Again, $\k_\zeta$ will be Galois over $\k$ and the prime
at $\infty$ will split completely.

\begin{questions}\label{roberts5}
Is there a canonical description of $\k_\zeta$ (perhaps as
multi-valued operators etc.) as a subfield of the
maximal totally split at $\infty$ subfield of $\k^{\rm sep}$?
\end{questions}

Next we examine the $v$-adic theory of multi-valued operators. Let $E$ be
a motive defined over a finite extension of $k_{v}$. One still has an
exponential $\exp_{v}$ and logarithm $\log_{v}$ associated to $E$. However,
it is important to note that one cannot expect that
$\exp_{v}$ is surjective.  In fact,
like the usual exponential function $p$-adically, the exponential of
the Carlitz module can easily be seen to 
have a finite $v$-adic radius of convergence. Still, there is still an obvious
notion of ``$v$-adic multi-valued operator'' on $M$ based on 
Definition \ref{multidef2} which we leave to the reader. 

To establish the $v$-adic analog of Theorem \ref{mainth} we first must
discuss some
results in the theory of hyper-derivatives.
Let $n$ be a positive integers and let $1\le j\le n$. Let $\Z_{\geq 0}$ be the
set of non-negative integers and let $P(n,j)$ be the set of those
$n$-tuples $(\mu_{1},\ldots,\mu_{n})\in (\Z_{\geq 0})^{n}$ such that
$ \sum_{i}\mu_{i}=j$ and $\sum_{i}i\mu_{i}=n$. For each $\mu\in P(n,j)$ one
defines an operator on $\A$ by
$${\mathfrak D}_{\mu}(f):={\mathfrak D}_{\mu_{1}}(f)^{\mu_{1}}\cdots
{\mathfrak D}_{\mu_{n}}(f)^{\mu_{n}}\,,$$
and, for $m\geq 0$, a multinomial coefficient 
$$M_{\mu}^{j}(m):=\frac{m!}{\mu_{1}!\cdots \mu_{n}! (m-j)!}=\frac{m(m-1)\cdots
(m-j+1)}{\mu_{1}!\cdots \mu_{n}!}\,.$$
We then have the following well-known result (see e.g., \cite{ky1}, \cite{us1})
which follows by induction.

\begin{lemma}\label{powerf}
Let $f\in \A$ and $m\geq 1$. Then
$${\mathfrak D}_{n}(f^{m})=\sum_{j=1}^{n}f^{m-j}\left(
\sum_{\mu\in P(n,j)}M_{\mu}^{j}(m){\mathfrak D}_{\mu}(f)\right)\,.$$
\end{lemma}

\begin{prop}\label{vadiccon}
Let $v\in {\rm Spec}(\A)$. Then the hyper-differential operators are
$v$-adically continuous.
\end{prop}

\begin{proof}
Let $n\geq 0$ and let ${\mathfrak D}_{n}$ be one such differential operator.
Let $v=(f)$ and let $g\in \A$ be $v$-adically small, say $g=cf^{m}$. We
need to establish that ${\mathfrak D}_{n}(g)$ is also $v$-adically
small. But Lemma \ref{powerf}
and the Leibniz formula easily show that ${\mathfrak D}_{n}(g)\equiv
0\pmod{f^{m-n}}$.
\end{proof}

Proposition \ref{vadiccon} immediately implies that the tangent action, 
$\alpha \mapsto E_{\alpha,\ast}$ for $\alpha\in \A$, continues $v$-adically to 
$\alpha\in \A_v$.
Inverting a non-zero non-unit element as in Proposition \ref{toK} then extends
this continuously to  all $\alpha\in \k_v$. Recall that $\theta$ gives
an isomorphism between $\k_v$ and $k_v$ where the prime $v$ is identified
with its isomorphic image. Thus the next result follows via the same 
arguments as in Proposition \ref{toKbar}.

\begin{theorem}\label{tokvbar}
The $T$-action of $E$ can be uniquely extended to an injection of
$\Fr$-algebras, $x\mapsto E_x$\,, of the separable closure $\k_v^{\rm sep}$
of $\k_v$ into the algebra of multi-valued operators. Moreover, if
$E_{x,\ast}$ is the induced action on the tangent space of the origin, then
$E_{x,\ast}=\theta(x)\cdot{\mathbf 1}+N_x=\overline{x}+N_{x}$ where $N_x$ is
nilpotent.\end{theorem}

We leave the appropriate version of Theorem \ref{tokvbar} in the case of
complex multiplication to the reader. Similarly, Theorem \ref{tokvbar}
leads to obvious $v$-adic versions of Questions \ref{bigbasic} and
\ref{multicompl}.

Theorem \ref{tokvbar} also raises some interesting questions about the
nature of multi-valued operators. We will pose these in the simplest
case of the tensor powers of the Carlitz module.
Let $C^{\otimes n}$, $n\geq 1$ be the $n$-th tensor power of the Carlitz
module viewed over $k$. By Proposition \ref{vadiccon} we know that we
may continuously extend $C^{\otimes n}_a$, $a\in \A$, to all 
$a\in \A_v$.
By our construction via the exponential and logarithm, it is clear that
all such $C^{\otimes n}_a$ belong to the formal algebra  of
formal power series $\sum a_i\tau^i$ where the $a_i$ are $n\times n$ matrices
with $k_v$-coefficients.
On the other hand, we can write $a$ as the limit of $a_i\in \A$ where
all the coefficients of $C^{\otimes n}_{a_i}$ are matrices with coefficients
in $A$. As such we immediately deduce that the coefficients of 
$C^{\otimes n}_a$ are matrices with $A_v$-valued coefficients for all
$a\in \A_v$. (See e.g., the discussion at the beginning of Section 4 of 
\cite{go1} for the case of Drinfeld modules.) 
We can thus {\it reduce} these elements modulo $v$ to get a formal
$\A_v$-action associated to $C^{\otimes n}$ over the {\it finite} field $A/v$. 
 
Now let $F$ be a finite separable extension of $\k$ and let
$\mathcal O$ be the $\A$-integers of $F$. Let $w$ be an unramified prime
of $\mathcal O$ of degree $1$ over  $v\in {\rm Spec}(\A)$; so ${\mathcal O}_w
\isom \A/v$.
Thus by the argument
just given there is a formal $\mathcal O$-action associated to
$C^{\otimes n}$ over
${\mathcal O}/w\isom \A/v$ by simply reducing modulo $w$.

\begin{questions} Is it possible to reduce the formal $\mathcal O$-action at
other primes?
\end{questions}

In other words, what ``denominators'' are induced by the construction of
multi-valued operators?

\end{document}